\documentclass[preprint,12pt]{elsarticle}

\usepackage{amsmath}
\usepackage{amssymb}
\usepackage{amsfonts}
\usepackage{amsbsy}
\usepackage{bm}
\usepackage{graphicx}
\usepackage{color} 
\usepackage[normalem]{ulem}

\usepackage{tensor}
\usepackage{tikz-cd}
\usepackage{stmaryrd}

\newcommand{\eg}{e.\,g.}%
\newcommand{\wrt}{w.r.t.}%

\newcommand{\R}{\mathbb{R}}%

\newcommand{\landau}{\mathcal{O}}

\newcommand{\bulk}{\mathcal{V}}

\newcommand{\surf}{\mathcal{S}}

\newcommand{\surfh}{\surf_{h}}

\newcommand{\charFct}[1]{\chi_{#1}}

\newcommand{\meshSize}{h}

\newcommand{\Torus}{\cal{T}}

\newcommand{\TangentBundle}[2]
{\mathrm{T}_{#1}{#2}}
\newcommand{\TensorBundle}[3][1={n},2={}]
{\mathrm{T}^{#1}_{#2}{#3}}

\newcommand{\normal}{\boldsymbol{\nu}}
\newcommand{\tangent}{\textup{T}}

\newcommand{\ext}[1]{\underline{#1}}

\newcommand{\fb}{\boldsymbol{f}}%
\newcommand{\xb}{\boldsymbol{x}}%
\newcommand{\Qb}{\boldsymbol{Q}}%
\newcommand{\ub}{\boldsymbol{u}}%
\newcommand{\Ub}{\boldsymbol{U}}%

\newcommand*{\rom}[1]{\textup{\uppercase\expandafter{\romannumeral#1}}}

\newcommand{\shapeOperator}{\boldsymbol{B}}%
\newcommand{\ProjSurf}{\operatorname{\Pi}}
\newcommand{\ProjClPoint}{\pi}

\newcommand{\gradSurf}{\mathrm{grad}\,}
\newcommand{\divSurf}{\mathrm{div}\,}

\newcommand{\LinLagrange}{\mathcal{L}^1}

\newcommand{\llScalIdx}[3]{\langle #1 , #2 \rangle_{#3}}
\newcommand{\llScalSurf}[2]{\llScalIdx{#1}{#2}{\surf}}
\newcommand{\llScalBulk}[2]{\llScalIdx{#1}{#2}{\bulk}}

\newcommand{\intBulkh}[1]{\int_{\bulk_{\meshSize}} #1 \, \mathrm{d}\bulk_{\meshSize}}

\newcommand{\sgnDist}{\rho}
\newcommand{\phase}{\phi}

\newcommand{\doubleWell}{W}

\newcommand{\epsPhase}{\epsilon}
\newcommand{\deltaPhase}{\delta}

\newcommand{\Velocity}{\Ub}
\newcommand{\velocity}{\ub}

\newcommand{\VeloTest}{\boldsymbol{\Psi}}

\newcommand{\nrmPenalty}{C_N}

\newcommand{\meshconv}{$\mathtt{meshconv}\;$}
\newcommand{\paraview}{$\mathtt{paraview}\;$}
\newcommand{\DIerrLtwo}{\mathrm{E}_{DI}}
\newcommand{\EXPLerrLtwo}{\mathrm{E}_{\surf}}

\begin{document}

\begin{frontmatter}
\title{A diffuse interface approach for vector-valued PDEs on surfaces}

\author[label1]{Michael Nestler}
\affiliation[label1]{organization={Institute of Scientific Computing, TU Dresden},
            city={Dresden},
            postcode={01062}, 
            country={Germany}}

\author[label1,label2,label3]{Axel Voigt}
\affiliation[label2]{organization={Center for Systems Biology Dresden (CSBD)},
            addressline={Pfotenhauerstra{\ss}e 108}, 
            city={Dresden},
            postcode={01307}, 
            country={Germany}}
\affiliation[label3]{organization={Cluster of Excellence Physics of Life (PoL)},
            city={Dresden},
            postcode={01062}, 
            country={Germany}}

\begin{abstract}
Approximating PDEs on surfaces by the diffuse interface approach allows us to use standard numerical tools to solve these problems. This makes it an attractive numerical approach. We extend this approach to vector-valued surface PDEs and explore their convergence properties. In contrast to the well-studied case of scalar-valued surface PDEs, the optimal order of convergence can only be achieved if certain relations between mesh size and interface width are fulfilled. This difference results from the increased coupling between the surface geometry and the PDE for vector-valued quantities defined on it.
\end{abstract}

\end{frontmatter}

\section{Introduction}
	
PDEs on surfaces remain an active field of research in applied mathematics and computational science. Due to their coupling with the geometry of such surfaces, PDEs are intrinsically nonlinear. This leads to new challenges in modeling and numerical analysis. Most of these challenges are addressed for scalar-valued surface PDEs, see \cite{dziuk2013finite} for a review. In the scalar case the coupling between surface geometry and the PDE is relatively weak, and thus numerical approaches established in flat spaces are applicable after small modifications. For vector-valued surface PDEs these approaches are no longer sufficient. Surface vector-fields often need to meet additional constraints. One example is the tangentiality of these fields. In this case, they need to be considered elements of the tangent bundle of the surface. This will lead to a strong nonlinear coupling between the surface geometry and the PDE. One break through which allows us to deal with these new challenges was the idea to express the solution and the surface differential operators in the global coordinate system of the embedding space and to penalize normal components, independently introduced in  \cite{nestler2018orientational,jankuhn2018,hansbo2020analysis} and generalized in \cite{nestler2019finite}. This approach essentially allows us to apply established tools for scalar-valued surface PDEs to each component. Popular approaches are surface finite elements (SFEM) and trace finite elements (TraceFEM), which have been applied to various problems in liquid crystal theory \cite{nestler2018orientational,nitschke2018nematic,PhysRevFluids.4.044002,nitschke2020liquid,nestler2020properties}, fluid mechanics \cite{doi:10.1063/1.5005142,jankuhn2018,ReutherEtAl2020FluidDeformableSurfaces,olshanskii2018finite,brandner2022finite} and biological physics \cite{CiCP-31-947,WMB_arXiv_2022,BKV_arXiv_2023}. Numerical analysis results for these methods also exist, but are restricted to the most simple equations of this type, the surface vector-valued Laplace equation \cite{grande2018analysis,hansbo2020analysis}, the surface vector-valued Helmholtz equation \cite{HarderingPraetorius2021TangentialErrors} and the surface Stokes equations \cite{OlshanskiiEtAl2019InfSupStability,JankuhnEtAl2021ErrorAnalysis}. All these results show the necessity for an appropriate approximation of the geometric quantities that enter these equations. Optimal order of convergence can often only be achieved if this approximation is of a higher order than the solution. These results reflect the increased coupling between the surface geometry and the solution for vector-valued surface PDEs if compared with scalar-valued surface PDEs. This increase in complexity is also shown numerically by comparing scalar-, vector- and tensor-valued surface diffusion equations by various numerical methods \cite{bachini2022diffusion}.

Another popular method to solve surface PDEs, the diffuse interface method \cite{ratz2006pde,burger2009finite} is less explored for vector-valued surface PDEs. Two examples, where solutions are compared with other more established methods are \cite{nestler2018orientational,bachini2022diffusion}. The diffuse interface method approximates the surface PDE by a bulk PDE that can be solved with standard numerical tools. This makes the approach attractive to be used in various application areas for more complex problems, especially those where the surface evolves according to physical laws which depend on the vector-valued surface quantity. Such applications can be found in biology, e.g. in morphogenesis, where the method recently led to spectacular results \cite{hoffmann2022theory}. However, knowing the subtleties associated with the approximation of geometric terms in SFEM and TraceFEM, discussed above, and the additional error emerging from the approximation of the diffuse interface method, such results should be considered with care. 

We here consider a tangential vector-valued surface Helmholtz equations to explore the convergence properties of the diffuse interface methods. This can be considered as a model problem, and the obtained results can be assumed to be applicable also to other more complex surface PDEs. In Section \ref{sec:model} we introduce the surface model, review the described extension to the embedding space and required penalization of the normal component, which provides the basis for a SFEM discretization and derive from this the diffuse interface formulation. We formulate these approaches in variational form and provide a finite element (FEM) discretization. The diffuse interface approximation is justified by formal matched asymptotics, which directly follow from the results for scalar-valued surface PDEs \cite{ratz2006pde}. In Section \ref{sec:results} we construct an analytical solution and perform various convergence studies. The results indicate the necessity for accurate approximations of the surface normals. These approximations can be obtained from the phase field variable or the signed distance function of the implicit surface description. Depending on the chosen approach specific minimal requirements regarding the relation of interface width and mesh size at the interface have to be meet for optimal convergence. Ignoring these requirements can result in numerical solution procedures that do not converge! These results significantly differ from the known results for the diffuse interface method for scalar-valued surface PDEs. In Section \ref{sec:conclusion} we draw conclusions.

\section{Model problem} \label{sec:model}
We demonstrate the applicability of the diffuse interface approach by considering the surface vector-valued Helmholtz equation on a smooth, oriented two-dimensional surface $\surf \subset \R^3$ without boundaries. The equation reads
\begin{align}
	-\divSurf \gradSurf \velocity + \velocity & = \fb \quad \mbox{on }\surf, \label{eq:StrVecHelmholtz}
\end{align}
for a tangential vector field $\velocity$ defined in the tangent bundle $\tangent\surf$ of $\surf$ and a compatible right hand side $\fb$. The considered differential operators $\gradSurf$ and $\divSurf$ correspond to the metric-preserving covariant derivative $\nabla_{\surf}$. The diffuse interface formulation will be based on the SFEM formulation for vector-valued surface PDE's, as detailed in \cite{nestler2019finite}. The general idea is to express the solution and the operators in the global coordinate system of the embedding space $\R^3$ and to penalize normal components. We denote the Euclidean coordinates of $\R^3$, by uppercase letters $I,J,K, \ldots$ and consider Einstein's sum convention. We consider the outward oriented normal $\normal$, the surface identity $\ProjSurf= \mathbb{I}-\normal\normal$ and the shape operator $\shapeOperator = -\nabla_{\surf} \normal$. For $\velocity$ we denote the formal extension to the embedding space by $\Velocity$, for which we require $[\Velocity - (\Velocity \cdot \normal) \normal ]|_{\surf}=\velocity$. Finally we define the full tangential projection of such extended fields along all their components by $\ProjSurf[.]$, \eg\ for a vector-valued field $\ProjSurf[{\mathbf{P}}]_{I} = \ProjSurf_{IJ} {P}_J$ and for a 2-tensor-valued field $\ProjSurf[{\Qb}]_{IJ} = \ProjSurf_{IK} {Q}_{KL} \ProjSurf_{LJ}$.  

Given these notations we can express the surface differential operators by their $\R^3$ counterparts. Acting on the extended field $\Velocity$ this reads
\begin{align}
	\gradSurf \velocity & = \ProjSurf[\nabla \Velocity]|_{\surf} + \shapeOperator(\Velocity|_{\surf}\cdot \normal) \\
	[\gradSurf \velocity]_{IJ} & = \ProjSurf_{IL} \partial_K \Velocity_L \ProjSurf_{KJ} + \shapeOperator_{IJ}\Velocity_{L}\normal_{L},
\end{align}
see \cite{nestler2019finite} for details. To obtain the embedded variational formulation we use scalar products of pointwise full contraction of the tangential parts of the extended solution $\Velocity$ and the test function $\VeloTest$, e.g.
\begin{align}
	\llScalSurf{\velocity}{\VeloTest} = \llScalSurf{\Velocity}{\VeloTest} = \int_{\surf} (\ProjSurf[\Velocity],\ProjSurf[\VeloTest])\,\mathrm{d}\surf = \int_{\surf} \Velocity_I \ProjSurf_{IJ}\VeloTest_J\,\mathrm{d}\surf.
\end{align}
From these ingredients we set up a component wise solution-test space $[H^1(\surf)]^3$ and write the embedded variational problem
\begin{align}
	\llScalSurf{\ProjSurf[\nabla \Velocity] + \shapeOperator(\Velocity\cdot \normal)}{\ProjSurf[\nabla \VeloTest] 
		+ \shapeOperator(\VeloTest\cdot \normal)}&\nonumber \\ + \llScalSurf{\Velocity}{\VeloTest} + \nrmPenalty\llScalSurf{\Velocity\cdot \normal}{\VeloTest \cdot \normal}& = \llScalSurf{\fb}{\VeloTest} &\forall \VeloTest \in [H^1(\surf)]^3 \label{eq:VarVecHelmholtz}
\end{align}
with an additional penalty term with prefactor $\nrmPenalty$ to approximate $\Velocity\cdot \normal = 0$. In this formulation each component can be considered by classical SFEM for scalar-valued fields, see \cite{dziuk2013finite}.

Various numerical studies confirm the applicability of this approach. Numerical analysis provides various estimates of solution error convergence for surface mesh size $\meshSize_{\surf} \rightarrow 0$. Here, we consider the surface $\surf$ to be approximated by faceted polyhedra $\surfh$. Such first order geometry approximation allows only linear Lagrangian elements $[\Velocity_{\meshSize}]_I \in \LinLagrange(\surfh) \approx H^1(\surf)$ for the discretization of the component solution-test space. This approach is referred to as the linear isogeometric ansatz. With the penalty prefactor $\nrmPenalty = \landau(1/\meshSize_{\surf}^2)$ one obtains quadratic convergence of the $L^2$-error $$\EXPLerrLtwo(\velocity -\ProjSurf[\Velocity_{\meshSize}];\meshSize_{\surf}) = \left(\llScalSurf{\velocity -\ProjSurf[\Velocity_{\meshSize}]}{\velocity - \ProjSurf[\Velocity_{\meshSize}]}\right)^{1/2}.$$ Please note, by approximating $\surf$ by a faceted polyhedra $\surfh$ we yield $\shapeOperator_{\meshSize} = 0$ which significantly simplifies the discretized weak formulation of the Laplacian in eq. \eqref{eq:VarVecHelmholtz}. Various numerical studies, also for more complicated problems \cite{nitschke2020liquid,CiCP-31-947,BKV_arXiv_2023}, confirm quadratic order.
	
At first glance it seems straight forward to apply the tools introduced in \cite{ratz2006pde} to each component in eq. \eqref{eq:VarVecHelmholtz} to obtain a diffuse interface approximation which turns the problem defined on $\surf$ into a problem defined in $\R^3$. However, eq. \eqref{eq:VarVecHelmholtz} contains various geometric terms, which are not considered in \cite{ratz2006pde} or any subsequent analysis. The sensitivity of the solution to the approximation of these geometric terms is known from \cite{hansbo2020analysis,HarderingPraetorius2021TangentialErrors} and contrary to the scalar case \cite{doi:10.1137/070708135}, higher order approximations for the surface and the solution only lead to better convergence properties if a higher order approximation of the normals is used in the introduced penalization term. So, the question arises if an additional approximation of the normals, resulting from an implicit description of the surface in the diffuse interface formulation, is sufficient to obtain the same convergence properties as the SFEM approach. We will answer this question in the following.
	
In order to formulate the diffuse interface approximation we extend the quantities $\normal$, $\ProjSurf$, $\shapeOperator$ and $\fb$ to the embedding space. We consider this component wise constant in the normal direction and denote the extended quantities by $\ext{\normal}$, $\ext{\ProjSurf}$, $\ext{\shapeOperator}$ and $\ext{\fb}$, respectively. With this we can write $\ext{\shapeOperator}_{IJ}=- \ext{\ProjSurf}_{IK} \partial_K \ext{\normal}_J$, and $ \ext{\shapeOperator}_{IJ}|_{\surf} = \shapeOperator_{IJ}$. We now follow \cite{ratz2006pde} and embed the surface in a bulk domain $\surf \subset \bulk$ and use an implicit description of $\surf$. For this purpose we use a signed distance function $\sgnDist(\xb)$ with $|\nabla \sgnDist | = 1$ and $\surf = \{\xb \in \bulk,\, \sgnDist(\xb)=0\}$ and a phase field function $\phase(\xb) = 1/2 \left(1 -\tanh(3\sgnDist(\xb)/\epsPhase)\right)$ with interface width $\epsPhase > 0$ and $\surf = \{\xb \in \bulk,\, \phase(\xb)=1/2\}$. With these functions we can evaluate extended normals. We obtain 
\begin{align}
	\mbox{[A]:}\; \ext{\normal} = \nabla \sgnDist, \quad  \mbox{ [B]:}\; \ext{\normal} = -\nabla \phase / | \nabla \phase \label{eq:nrmApprox}
\end{align}
The surface delta function is approximated by $\charFct{\surf} \approx \doubleWell(\phase) = 36/\epsPhase \,(1-\phase)^2 \phase^2$, the classical double-well potential in Ginzburg-Landau energies \cite{cahn1958free}. Within this notation we substantiate the previously used constant normal extension by the closest point projection $\ProjClPoint: \bulk \rightarrow \surf, \, \ProjClPoint(\xb) = \xb - \sgnDist \ext{\normal}$, \eg\ $\ext{\fb}(\xb) = \fb \circ \ProjClPoint(\xb)$, which is well defined for $\epsPhase$ sufficient small such that $\epsPhase \|\shapeOperator(\xb)\| < 1, \forall \xb \in \surf$. With this notation we express the surface scalar product by a scalar product of the embedding space
\begin{align}
	\llScalSurf{\velocity}{\VeloTest} = \llScalSurf{\Velocity}{\VeloTest} \approx \int_{\bulk} \doubleWell(\phase) \Velocity_I \ext{\ProjSurf}_{IJ}\VeloTest_J\,\mathrm{d}\bulk = \llScalBulk{\Velocity}{\VeloTest}.
\end{align}
We use this embedded scalar product to set up the embedded variational problem with a component wise solution space $[H^1(\bulk)]^3$. So the diffuse interface approximation of the component wise embedded variational problem eq. \eqref{eq:VarVecHelmholtz} reads
\begin{align}
	\llScalBulk{\left[\ext{\ProjSurf}[\nabla \Velocity] + \ext{\shapeOperator}(\Velocity\cdot \ext{\normal})\right]}{\ext{\ProjSurf}[\nabla \VeloTest] + \ext{\shapeOperator}(\VeloTest\cdot \ext{\normal})}& \\
	 + \llScalBulk{\Velocity}{\VeloTest} + \nrmPenalty\llScalBulk{\left[\Velocity\cdot \ext{\normal}\right]}{\left[\VeloTest \cdot \ext{\normal}\right]} 	&= \llScalBulk{\ext{\fb}}{\VeloTest} \quad \forall \VeloTest \in [H^1(\bulk)]^3. \nonumber \label{eq:DIVarVecHelmholtz}
\end{align}
Considering such component wise formulation as a coupled system of scalar fields, the matched asymptotic analysis of \cite{ratz2006pde} provides formal convergence of eq. \eqref{eq:DIVarVecHelmholtz} to eq. \eqref{eq:VarVecHelmholtz} for $\epsPhase \to 0$. This establishes eq. \eqref{eq:DIVarVecHelmholtz} to be a diffuse interface approximation of eq. \eqref{eq:StrVecHelmholtz}. Using this procedure any surface PDE can be reformulated into the corresponding diffuse interface approximation.
	
In analogy to the SFEM approach, we use a linear ansatz for the implicit geometric description by $\sgnDist_{\meshSize}, \phase_{\meshSize} \in \LinLagrange(\bulk_{\meshSize})$ for a tetrahedral mesh $\bulk_{\meshSize}$ of $\bulk$, where $\meshSize$ denotes the bulk mesh size in the interface $\phase \in [0.05, 0.95]$. Within such ansatz the recovered surface $\surfh = \{ \sgnDist_{\meshSize} \equiv 0\}$ or $\surfh = \{ \phase_{\meshSize} \equiv 1/2\}$ is a faceted polyhedra. The same ansatz is considered for the solution $\Velocity_{\meshSize} \in [\LinLagrange(\bulk_{\meshSize})]^3$. This allows us to consider $\ext{\shapeOperator}_{\meshSize} = 0$ across the interface and we obtain 
\begin{equation}
	H(\Velocity_{\meshSize},\VeloTest) + P(\Velocity_{\meshSize},\VeloTest) + R(\Velocity_{\meshSize},\VeloTest) = F(\VeloTest) \quad \forall\,  \VeloTest \in [\LinLagrange(\bulk_h)]^3 \label{eq:DI_VecHelmholtz}
\end{equation}
with the Helmholtz operator $H(\Velocity_{\meshSize},\VeloTest)$, the normal penalty term $P(\Velocity_{\meshSize},\VeloTest)$, a bulk stabilization $R(\Velocity_{\meshSize},\VeloTest)$ and the right hand side $F(\VeloTest)$, defined by:
\begin{align}
	H(\Velocity_{\meshSize},\VeloTest) =&  \intBulkh{\doubleWell(\phase_{\meshSize})[\ext{\ProjSurf}_{\meshSize}]_{IJ} \; \partial_K \Velocity_{\meshSize,I}\partial_K \VeloTest_J} \nonumber \\
	& + \intBulkh{\doubleWell(\phase_{\meshSize})[\ext{\ProjSurf}_{\meshSize}]_{IJ}\;\Velocity_{\meshSize,I}\VeloTest_J} \nonumber \\
	P(\Velocity_{\meshSize},\VeloTest) = & \, \nrmPenalty \intBulkh{ \doubleWell(\phase_{\meshSize}) [\ext{\normal}_{\meshSize}]_I \Velocity_{\meshSize,I}[\ext{\normal}_{\meshSize}]_{J} \VeloTest_J} \nonumber \\
	R(\Velocity_{\meshSize},\VeloTest) = & \, \deltaPhase \intBulkh{\partial_K \Velocity_{\meshSize,I}\partial_K \VeloTest_I} \nonumber \\
	F(\VeloTest) = & \intBulkh{\doubleWell(\phase_{\meshSize})[\ext{\ProjSurf}_{\meshSize}]_{IJ}\;\ext{\fb}_I\VeloTest_J}, \nonumber
\end{align}
respectively. Here we have used a component wise constant regularization outside the interface of $\deltaPhase=10^{-6}$ as discussed for the scalar case in \cite{ratz2006pde}. Inspired by the normal penalty condition, specified in \cite{HarderingPraetorius2021TangentialErrors} for the isogeometic case, we consider also within the diffuse interface approximation $\nrmPenalty = 10/h^2$.
	

	
\section{Results} \label{sec:results}
	
We are interested in the convergence properties of the diffuse interface approach. We have to take into consideration that error measures will depend on $\epsPhase$ and $\meshSize$. For further discussion we define a relation $(\meshSize, \epsPhase)$ to essentially express the number of mesh points across the interface for a given interface width. For our considerations we follow \cite{feng2003numerical} and group those relations as linear, where $\meshSize = \landau(\epsPhase)$ so the number of points remains constant within the interface for $\epsPhase \rightarrow 0$, and higher order relations, where $\meshSize^2 = \landau(\epsPhase^3),\landau(\epsPhase^4),\landau(\epsPhase^5)\hdots$ such that the number of points within the interface increases for $\epsPhase \rightarrow 0$, see Figure-\ref{fig:TorusMeshModes}-[B].
	
Within the same setting, a linear ansatz for the geometric description and a linear ansatz for the solution, in \cite{ratz2006pde} a numerical study was performed to estimate convergence rates. For a linear relation $(\meshSize, \epsPhase)$ quadratic convergence was obtained for the scalar-valued problem, which can be considered optimal. It thus provides an upper limit to our vector-valued problem. However, it remains open if this limit can be reached in typical applications, where the signed distance function $\sgnDist_h$ or the phase field variable $\phase_h$ need to be constructed from a surface mesh or where the phase field variable $\phase_h$ might be a solution to another PDE. In these cases the computation of the normals $\normal_h$ by eq. \eqref{eq:nrmApprox}[A] or [B] adds an additional source of error. 
	
\subsection{Benchmark formulation}
\begin{figure}[ht!]
		\begin{center}
			\includegraphics[width=0.8\linewidth]{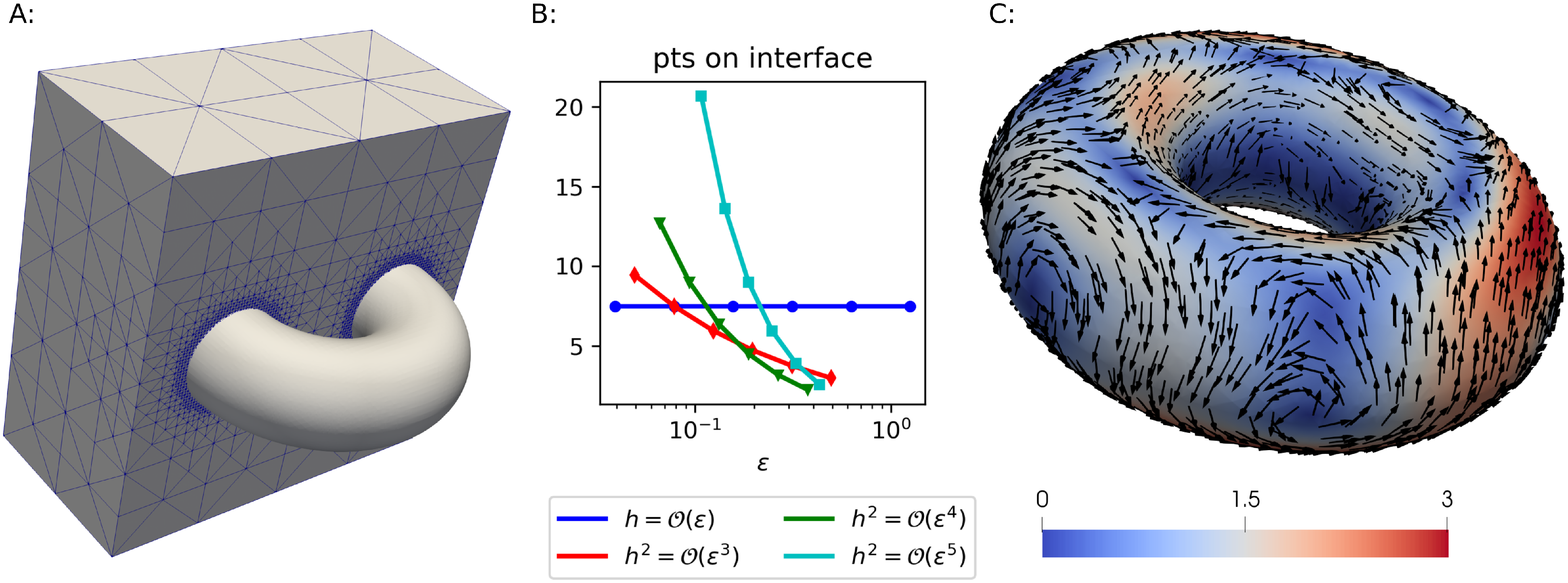}
		\end{center}
		\caption{\textbf{Adaptive meshes for torus embedded in box domain:}[A]: Meshed box domain $\bulk_{\meshSize}$ clipped at $x\equiv 0$ to visualize adaptive mesh with $\meshSize=0.03125$, plotted torus recovered by $\sgnDist_{\meshSize} \equiv 0$. [B]: Development of number mesh points across the interface $\phase \in [0.05,0.95]$ along $\epsPhase \rightarrow 0$ depending on $(\meshSize,\epsPhase)$ relation. [C]:Solenoidal vector field $\velocity$ on Torus, color encodes magnitude and glyphes direction of $\velocity$.}
		\label{fig:TorusMeshModes}
\end{figure}

To assess the impact of the major sources of approximations in the diffuse interface approach we perform a series of numerical experiments on a torus $\Torus$ with radi $R=1$ and $r=0.5$. The torus is embedded in $\bulk = [-2,2]^3$. The analytical signed distance function for this geometry is given in $\mathbb{R}^3,\, \xb = [x,y,z]$ by
\begin{align}
	\sgnDist(\xb) = \left(\left( \sqrt{x^2 + y^2}  - R \right)^2 + z^2 \right)^{1/2} - r \label{eq:sgnDistTorus}
\end{align}
The analytic phase field function $\phase$ is defined as before and for the analytic normal to $\Torus$ we consider $\normal = \nabla \sgnDist$. On $\Torus$ we consider the solenoidal tangential vector field $\velocity = \normal \times \nabla_{\surf}(x^2 y - 5 z^3)$ and its extension to $\bulk$ by $\Velocity = \velocity \circ \ProjClPoint$, see Figure-\ref{fig:TorusMeshModes}-[C]. This vector field will be used as an analytical solution for estimating the rate of error convergence. A compatible right hand side in eq. \eqref{eq:DI_VecHelmholtz} is constructed by
\begin{align}
	    F(\VeloTest) = \intBulkh{\doubleWell(\phase) \ext{\left[ \gradSurf \velocity \right] }_{IJ} \ext{\ProjSurf}[\nabla \VeloTest]_{IJ}} + \intBulkh{\doubleWell(\phase)\Velocity_I\ext{\ProjSurf}[\VeloTest]_I}.
\end{align}
To obtain results, which are comparable to the SFEM results in \cite{HarderingPraetorius2021TangentialErrors} we estimate $\| \Velocity - \ProjSurf[\Velocity_{\meshSize}]\| \leq \| \Velocity - \Velocity_{\meshSize}\| + \| \Velocity_{\meshSize} \cdot \normal \|$ and approximate surface error measured by the following diffuse interface $L^2$ error measures.
\begin{align}
	\DIerrLtwo(\Velocity - \Velocity_h; \meshSize, \epsPhase)& = \left(\intBulkh{\doubleWell(\phase) (\Velocity - \Velocity_h)\cdot(\Velocity - \Velocity_h)} \right)^{1/2} \\
	\DIerrLtwo(\Velocity_h \cdot \normal; \meshSize, \epsPhase)& = \left(\intBulkh{\doubleWell(\phase) \left[\Velocity_h\cdot \normal\right]\left[\Velocity_h \cdot \normal\right]} \right)^{1/2}.
\end{align}
For the subsequent numerical experiments we consider a set of mesh sizes $\meshSize \in [(1/2)^7, (1/2)^2]$, corresponding to bisections of $\bulk$. For each $\meshSize$ and $(\meshSize,\epsPhase)$ relation we evaluate an associated $\epsPhase$ and construct an adaptive tetrahedral mesh with meshsize $\meshSize$ in the diffuse interface of width $\epsPhase$. Off the interface, the mesh is only refined to preserve conformity. See Figure \ref{fig:TorusMeshModes}-[A]. We use the \meshconv toolbox \cite{meshconv} to generate these meshes and evaluate $\sgnDist_\meshSize$. \meshconv uses a triangulated surface as input which we generate via \paraview\cite{ahrens200536, ayachit2015paraview} with vertices on $\Torus$ and meshsize $\meshSize_{\surf} = \meshSize/4$. Different methods exist to compute signed distance functions, see e.g. \cite{Sussmanetal_JCP_1994, Russoetal_JCP_2000,Bornemannetal_CVS_2006,lee2017revisiting,royston2018parallel}. In \meshconv a second order method is implemented. We use the finite element toolbox \texttt{AMDiS} \cite{vey2007amdis,witkowski2015software} to solve eq. \eqref{eq:DI_VecHelmholtz}, which uses an unpreconditioned BiCGstabL solver provided by \texttt{PETSc} \cite{petsc-user-ref} for the linear system.
	
\subsection{Numeric evaluation of diffuse interface variables and normals} 
\begin{figure}[ht!]
		\begin{center}
			\includegraphics[width=\linewidth]{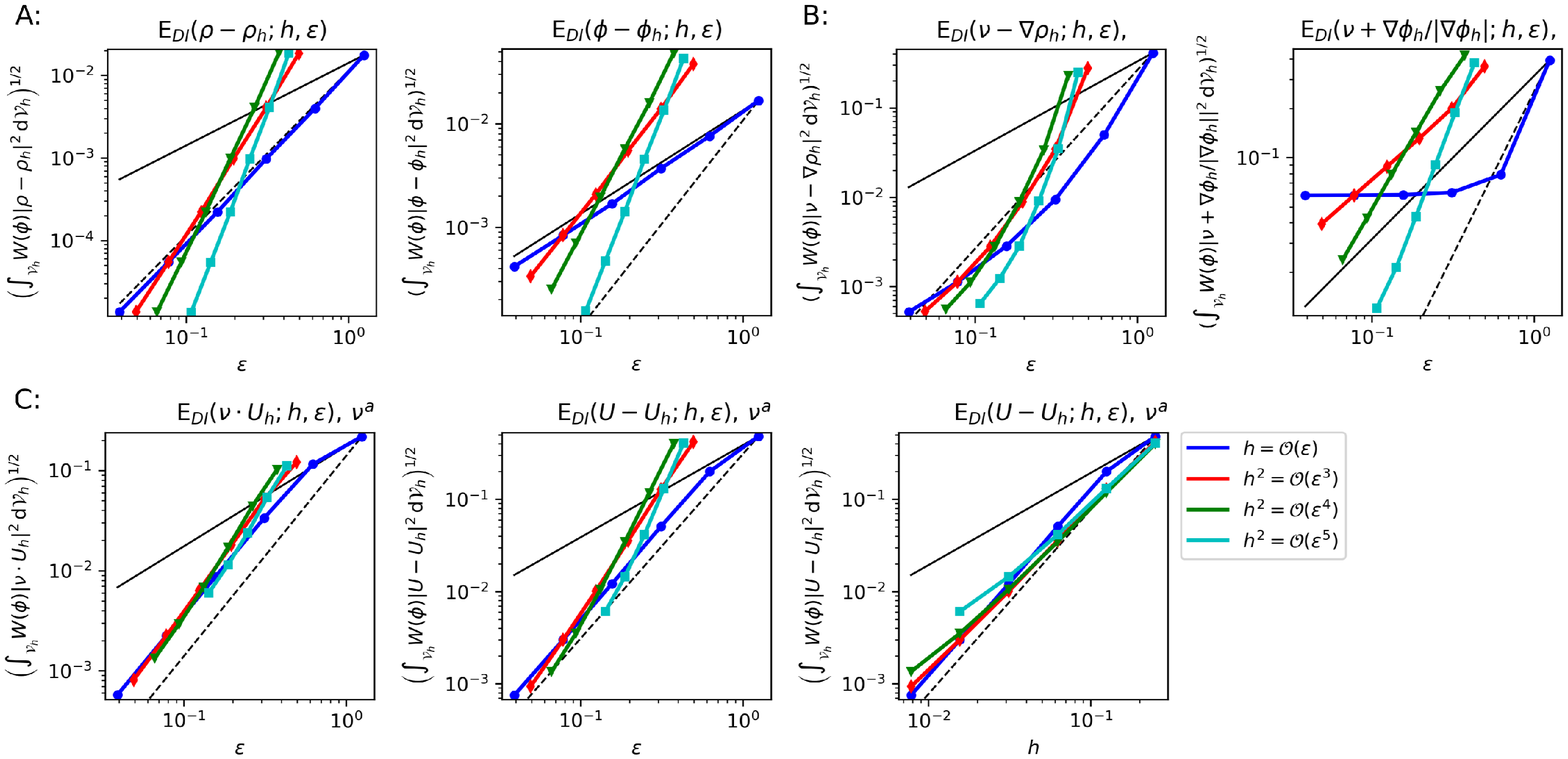}
		\end{center}
		\caption{\textbf{Error convergence for diffuse interface variables $\sgnDist_h$, $\phase_h$ and normal approximations $\normal_h$ and for solution $\Velocity_h$ with analytic diffuse interface variables:} [A]: Convergence of errors in evaluated diffuse interface variables, provided by \meshconv with $\epsPhase \rightarrow 0$, for signed distance $\sgnDist_h$ and phase field $\phase_h$ on $\bulk_h$. [B]: Convergence of errors in numeric approximation of normals depending on approximation method. (left) basing on signed distance with $\normal_h = \nabla \sgnDist_{\meshSize}$ and (right) derived from phase field $\normal_h = -\nabla \phase_{\meshSize} / | \nabla \phase_{\meshSize} |$ .[C]: Convergence \wrt interface width $\epsPhase$ for several $(\meshSize,\epsPhase)$ relations in solution normal component $\Velocity_h\cdot \normal$ (left), solution error for solving eq. \eqref{eq:DI_VecHelmholtz} with analytic $\phase$ and $\normal$ (middle) and solution errors versus mesh size $\meshSize$ (right).}
		\label{fig:TorusProblemSetup}
\end{figure}

First we examine the convergence of errors in the chosen method to compute diffuse interface variables $\sgnDist_\meshSize,\, \phase_\meshSize$ and $\normal_\meshSize$ derived from these fields. We compute $\sgnDist_\meshSize$ and $\phase_\meshSize$ with \meshconv for all $\meshSize$ and $(\meshSize,\epsPhase)$ relations. We evaluate the diffuse interface $L^2$-errors \wrt\ the analytical results $\DIerrLtwo(\sgnDist - \sgnDist_\meshSize; \meshSize, \epsPhase)$, and $\DIerrLtwo(\phase - \phase_\meshSize; \meshSize, \epsPhase)$ and plot these quantities versus $\epsPhase$, see Figure-\ref{fig:TorusProblemSetup}-[A]. We observe the errors in evaluating $\sgnDist_h$ converge quadratic, \wrt\ $\epsPhase \rightarrow 0$ for linear $(\meshSize,\epsPhase)$ relation. Higher order $(\meshSize,\epsPhase)$ relations increase the rate of convergence. For $\phase_h$ we observe the rates of convergence to be reduced about one order compared to rates for $\sgnDist_h$. Recalling that $|\nabla \sgnDist | = 1$ while $|\nabla \phase | \approx 1/\epsPhase$ near the interface such reduced order of convergence is expected. Using the evaluated diffuse interface fields to determine the normals along the approaches described in eq. \eqref{eq:nrmApprox} we evaluate the diffuse interface $L^2$-errors, $\DIerrLtwo(\normal - \nabla\sgnDist_h; \meshSize, \epsPhase)$ and $\DIerrLtwo(\normal + \nabla\phase_h / | \nabla \phase_h |; \meshSize, \epsPhase)$, see Figure-\ref{fig:TorusProblemSetup}-[B]. Using the linear $(\meshSize,\epsPhase)$ relation, see Figure-\ref{fig:TorusProblemSetup}-[B](blue lines), we yield quadratic rate of convergence in $\DIerrLtwo(\normal -\nabla \sgnDist_{\meshSize}; \meshSize, \epsPhase)$ while $\DIerrLtwo(\normal^a + \nabla \phase_{\meshSize} / | \nabla \phase_{\meshSize} |; \meshSize, \epsPhase)$ does not converge for $\epsPhase \rightarrow 0$. Considering higher order $(\meshSize,\epsPhase)$ relations the convergence rates, \wrt\ $\epsPhase$, for $\DIerrLtwo(\normal - \nabla\sgnDist_{\meshSize}; \meshSize, \epsPhase)$ could be improved beyond quadratic rate. In the case of $\DIerrLtwo(\normal +\nabla \phase_{\meshSize} / | \nabla \phase_{\meshSize} |; \meshSize, \epsPhase)$ a linear convergence rate could be achieved for $\meshSize^2 = \landau(\epsPhase^3)$ and a quadratic rate is achieved with $\meshSize^2 = \landau(\epsPhase^5)$ relation. We conclude that, in order to obtain normals with converging errors for $\epsPhase \rightarrow 0$, \meshconv provides sufficient good approximations for $\sgnDist_h$ to evaluate $\normal_h = \nabla \sgnDist_h$ even for a linear $(\meshSize,\epsPhase)$ relation. Yet considering the approach of computing normals from $\phase_h$, we have to keep in mind that the rate of convergence, \wrt\ $\epsPhase \rightarrow 0$, is effectively reduced by two orders. Therefore, errors in normals only converge with higher order $(\meshSize,\epsPhase)$ relations.
	
\subsection{Error convergence for solution with analytic diffuse interface variables} Next we consider the analytical forms for $\sgnDist$, $\phase$ and $\normal$ and solve eq. \eqref{eq:DI_VecHelmholtz} for a set of $\meshSize$ and $(\meshSize,\epsPhase)$ relations. We evaluate the diffuse interface $L^2$-error $\DIerrLtwo(\Velocity - \Velocity_h; \meshSize, \epsPhase)$, see Figure-\ref{fig:TorusProblemSetup}-[C]. We observe that for a linear $(\meshSize,\epsPhase)$ relation a quadratic rate of error convergence \wrt\ $\epsPhase \rightarrow 0$. This confirms the result from \cite{ratz2006pde} for scalar-valued PDEs. Again, using higher order $(\meshSize,\epsPhase)$ relations improves the rate of error convergence \wrt\ $\epsPhase \rightarrow 0$. But reviewing the error convergence \wrt\ $\meshSize$ it becomes apparent that $\epsPhase$ is the crucial factor in convergence of $\DIerrLtwo(\Velocity - \Velocity_h; \meshSize, \epsPhase)$ and using higher order $(\meshSize,\epsPhase)$ relations does not improve the rate of convergence for $\meshSize \rightarrow 0$. So, we conclude that an analytical diffuse interface approximation does not impair the convergence of errors \wrt\ $\epsPhase \rightarrow 0$, and we observe a quadratic rate of convergence for linear $(\meshSize,\epsPhase)$ relation.
	
\subsection{Error convergence for diffuse interface approximation} 
\begin{figure}[ht!]
		\begin{center}
			\includegraphics[width=0.8\linewidth]{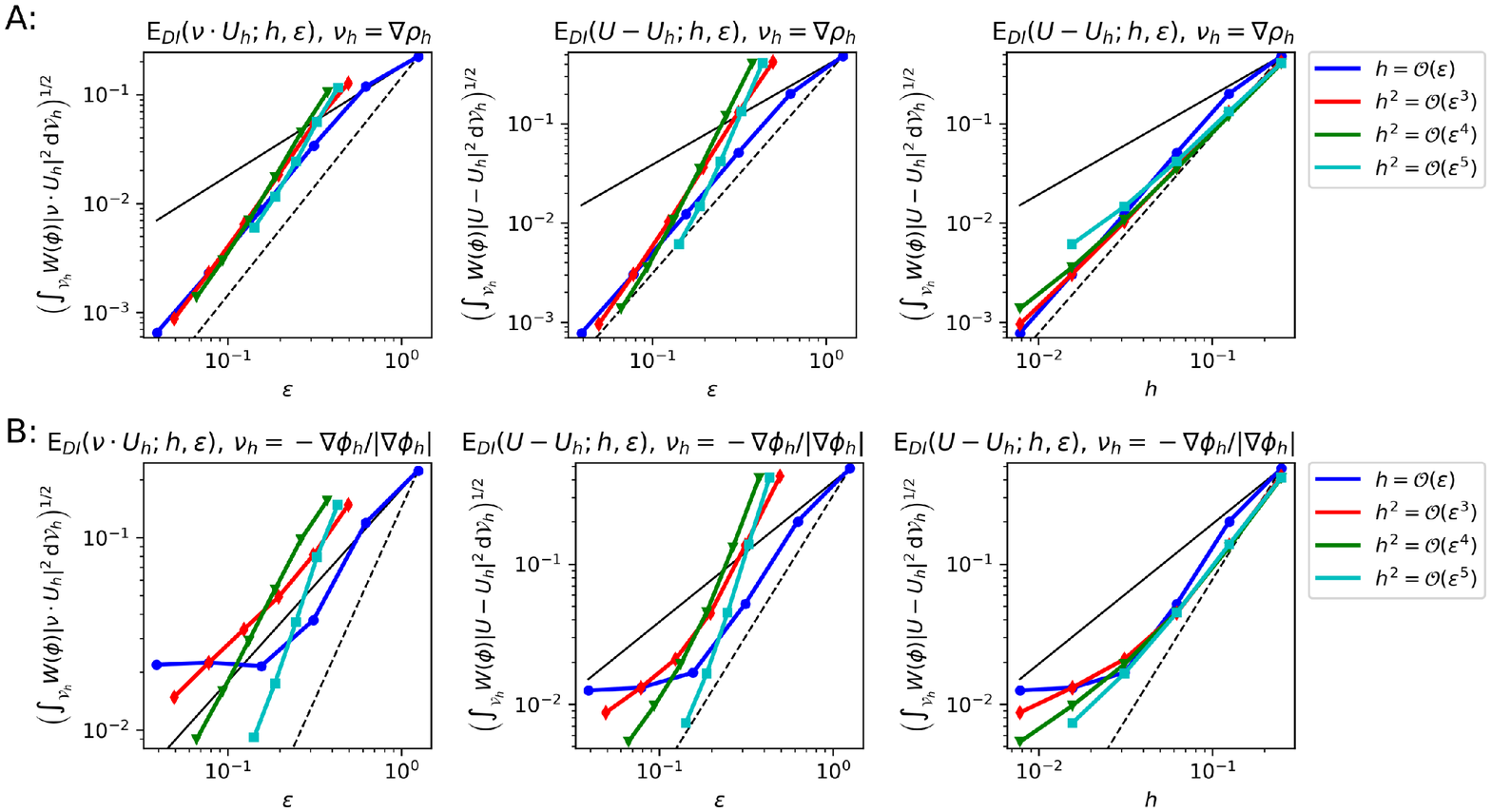}
		\end{center}
		\caption{\textbf{Error convergence for diffuse interface method for tangential surface vector-valued Helmholtz equation on torus:} [A]: Numerical experiments for $\epsPhase \rightarrow 0$ for normal approximation $\normal_h = \nabla \sgnDist_h$. (left) Error in incomplete suppression of normal components in solution $\Velocity_h$. (middle) Error convergence of $\DIerrLtwo(\Velocity -\Velocity_h; \meshSize, \epsPhase)$ \wrt\ $\epsPhase \rightarrow 0$. (right) Error convergence of $\DIerrLtwo(\Velocity -\Velocity_h; \meshSize, \epsPhase)$ \wrt\ $\meshSize \rightarrow 0$. [B]: Numerical experiments for $\epsPhase \rightarrow 0$ for normal approximation $\normal_h = -\nabla \phase_h/|\nabla \phase_h|$. (left) Error in incomplete suppression of normal components in solution $\Velocity_h$. (middle) Error convergence of $\DIerrLtwo(\Velocity -\Velocity_h; \meshSize, \epsPhase)$ \wrt\ $\meshSize \rightarrow 0$. (right) Error convergence of $\DIerrLtwo(\Velocity -\Velocity_h; \meshSize, \epsPhase)$ \wrt\ $\meshSize \rightarrow 0$. }
		\label{fig:TorusVectorHelmholtz}
\end{figure}

We now combine both, consider approximations of the diffuse interface variables $\sgnDist_\meshSize,\, \phase_\meshSize$ and $\normal_\meshSize$ obtained with \meshconv and solve eq. \eqref{eq:DI_VecHelmholtz} for a set of $\meshSize$ and $(\meshSize,\epsPhase)$ relations. For these numeric solutions $\Velocity_h$ we evaluate the diffuse interface $L^2$-error  $\DIerrLtwo(\Velocity - \Velocity_h; \meshSize, \epsPhase)$ and the normal component $\DIerrLtwo(\Velocity_h \cdot \normal; \meshSize, \epsPhase)$, see Figure-\ref{fig:TorusVectorHelmholtz}. For the normal approximation by $\normal_h = \nabla \sgnDist_h$, see Figure \ref{fig:TorusVectorHelmholtz}-[A] and a linear $(\meshSize,\epsPhase)$ relation we observe quadratic rate of convergence for $\epsPhase \rightarrow 0$ in the normal component and for the solution error. Due to the linear nature of the $(\meshSize,\epsPhase)$ relation this behavior translates into quadratic rate of convergence \wrt\ $\meshSize$ with absolute values very close to the results with analytic diffuse interface variables and normals in Figure \ref{fig:TorusProblemSetup}-[C]. For higher order $(\meshSize,\epsPhase)$ relations the errors and convergence rates are also very close to the results with analytic diffuse interface variables and normals. Quite contrary, we observe while using the approximation of $\normal_h = - \nabla \phase_h / | \nabla \phase_h |$ the errors in normal approximation to be dominant. Reviewing Figure \ref{fig:TorusProblemSetup}-[B](right) we observe that the reduced normal approximation quality of this approach yields reduced convergence rate for controlling the solutions normal component, see Figure \ref{fig:TorusVectorHelmholtz}-[B](left) directly impacting the solution error convergence rate in Figure \ref{fig:TorusVectorHelmholtz}-[B](middle). In the case of a linear $(\meshSize,\epsPhase)$ relation $\meshSize=\landau(\epsPhase)$ we observe no convergence in solution error, at $\meshSize^2=\landau(\epsPhase^3)$ we obtain a linear rate of convergence and for  $\meshSize^2=\landau(\epsPhase^4)$ and $\meshSize^2=\landau(\epsPhase^5)$ we yield close to quadratic rates of convergence. Comparing to the results for analytic diffuse interface variables and normals we observe only for $\meshSize^2=\landau(\epsPhase^5)$ similar errors and rates of convergence for $\DIerrLtwo(\Velocity_h \cdot \normal ; \meshSize, \epsPhase)$ and $\DIerrLtwo(\Velocity - \Velocity_h; \meshSize, \epsPhase)$ \wrt\ $\epsPhase \rightarrow 0$ and $\meshSize \rightarrow 0$. Therefore, we conclude that the diffuse interface method, based on component-wise FEM, yields rates of error convergence compatible to the results of SFEM for vector-valued surface PDEs only if a high-quality signed-distance function is used to compute the normals. In this case a linear $(\meshSize,\epsPhase)$ relation is sufficient. Alternatively if the normals are approximated by $\normal_h = - \nabla \phase_h / | \nabla \phase_h |$ a $(\meshSize,\epsPhase)$ relation with $\meshSize^2=\landau(\epsPhase^5)$ is required to achieve quadratic rate of convergence. A linear $(\meshSize,\epsPhase)$ relation does not lead to convergence in this case! This significantly differs from the diffuse interface method for scalar-valued surface PDEs introduced in \cite{ratz2006pde}.
	
\section{Discussion and conclusion} \label{sec:conclusion}
We provided a systematic derivation of a diffuse interface method to solve vector-valued surface PDEs. This approach builds on established SFEM formulations \cite{nestler2019finite} and allows us to consider each component as a scalar quantity. The resulting system of equations can be solved by standard FEM. The connection to the surface PDE is established by formal matched asymptotic following \cite{ratz2006pde} and is applicable to any vector-valued surface PDE. Along the same lines the approach can also be extended to tensor-valued surface PDEs.

Our numerical studies show that the same convergence properties as for SFEM can be achieved. With a linear isogeometric approach, which in our approach corresponds to linear Lagrange elements and a simplicial mesh, a quadratic order of convergence can be obtained for the considered vector-valued surface Helmholtz equation on a torus. However, to achieve this optimal order requires an accurate approximation of the surface normals. If they have to be numerically constructed a specific $(\meshSize,\epsPhase)$ relation is necessary. If the computation of the normals is based on the signed distance function $\sgnDist_h$ a linear relation is sufficient. However, if the computation of the normals is based on the phase field variable $\phase_h$ a linear relation does not lead to convergence and a higher order relation is required. The numerical results suggest $\meshSize^2=\landau(\epsPhase^5)$ to obtain a quadratic rate of convergence.
    
This result has severe consequences for applications. While for stationary surfaces, the computation of a signed distance function is not an issue and in such cases the diffuse interface method with normals computed from the signed distance function provides a valuable tool with optimal order of convergence, for moving surfaces this is no longer an option. The increased computational cost resulting from the computation of a signed distance function in each time step makes the approach inefficient. Assuming that the evolution of the surface is governed by a phase field problem which determines $\phase_h$ in each time step, the normals should be computed using $\phase_h$. Such approaches are common for scalar-valued surface quantities, such as concentration \cite{burger2006surface,ratz2006diffuse,lowengrub2009phase,teigen2011diffuse,marth2014signaling,doi:10.1137/21M1433642} or particle density \cite{aland2012buckling,aland2011continuum,PhysRevE.86.046321}. However, for vector-valued surface quantities using this approach does not lead to convergence if a linear $(\meshSize,\epsPhase)$ relation is considered. We can only expect convergence for higher order relations and quadratic convergence for $\meshSize^2=\landau(\epsPhase^5)$. This needs to be respected in applications. Accounting for this issue increases the computational cost. However, the advantages of using basic numerical tools like linear Lagrangian elements and generic iterative solvers and the possibility to deal with topological changes in the evolution of the surface remain. This together with the possibility to couple the surface PDE with equations in the embedding space \cite{li2009solving,teigen2009diffuse}, e.g. in two-phase flow problems with fluidic surfaces \cite{reuther2016incompressible}, makes this approach appealing for various applications.\\
    
{\bf{Acknowledgments:}} 
This work was supported by the German Research Foundation (DFG) through grant VO 899/22 within the Research Unit ”Vector- and Tensor-Valued Surface PDEs” (FOR 3013). We further acknowledge computing resources provided by ZIH at TU Dresden and by JSC at FZ J\"ulich, within projects WIR and PFAMDIS, respectively.
	
	
	\bibliographystyle{elsarticle-num}
	\bibliography{lib} 

\begin{thebibliography}{10}
\expandafter\ifx\csname url\endcsname\relax
  \def\url#1{\texttt{#1}}\fi
\expandafter\ifx\csname urlprefix\endcsname\relax\def\urlprefix{URL }\fi
\expandafter\ifx\csname href\endcsname\relax
  \def\href#1#2{#2} \def\path#1{#1}\fi

\bibitem{dziuk2013finite}
G.~Dziuk, C.~M. Elliott, Finite element methods for surface {PDEs}, Acta
  Numerica 22 (2013) 289--396.

\bibitem{nestler2018orientational}
M.~Nestler, I.~Nitschke, S.~Praetorius, A.~Voigt, Orientational order on
  surfaces: The coupling of topology, geometry, and dynamics, Journal of
  Nonlinear Science 28 (2018) 147--191.

\bibitem{jankuhn2018}
T.~Jankuhn, M.~A. Olshanskii, A.~Reusken, {Incompressible fluid problems on
  embedded surfaces: Modeling and variational formulations}, Interfaces and
  Free Boundaries 20 (2018) 353--377.
\newblock \href {https://doi.org/{10.4171/IFB/405}}
  {\path{doi:{10.4171/IFB/405}}}.

\bibitem{hansbo2020analysis}
P.~Hansbo, M.~G. Larson, K.~Larsson, Analysis of finite element methods for
  vector laplacians on surfaces, IMA Journal of Numerical Analysis 40 (2020)
  1652--1701.

\bibitem{nestler2019finite}
M.~Nestler, I.~Nitschke, A.~Voigt, A finite element approach for vector-and
  tensor-valued surface pdes, Journal of Computational Physics 389 (2019)
  48--61.

\bibitem{nitschke2018nematic}
I.~Nitschke, M.~Nestler, S.~Praetorius, H.~L{\"o}wen, A.~Voigt, Nematic liquid
  crystals on curved surfaces: a thin film limit, Proceedings of the Royal
  Society A 474~(2214) (2018) 20170686.

\bibitem{PhysRevFluids.4.044002}
I.~Nitschke, S.~Reuther, A.~Voigt, Hydrodynamic interactions in polar liquid
  crystals on evolving surfaces, Phys. Rev. Fluids 4 (2019) 044002.
\newblock \href {https://doi.org/10.1103/PhysRevFluids.4.044002}
  {\path{doi:10.1103/PhysRevFluids.4.044002}}.

\bibitem{nitschke2020liquid}
I.~Nitschke, S.~Reuther, A.~Voigt, Liquid crystals on deformable surfaces,
  Proceedings of the Royal Society A 476 (2020) 20200313.

\bibitem{nestler2020properties}
M.~Nestler, I.~Nitschke, H.~L{\"o}wen, A.~Voigt, Properties of surface
  {Landau--de Gennes Q}-tensor models, Soft Matter 16 (2020) 4032--4042.

\bibitem{doi:10.1063/1.5005142}
S.~Reuther, A.~Voigt, Solving the incompressible surface {Navier-Stokes}
  equation by surface finite elements, Physics of Fluids 30 (2018) 012107.
\newblock \href {https://doi.org/10.1063/1.5005142}
  {\path{doi:10.1063/1.5005142}}.

\bibitem{ReutherEtAl2020FluidDeformableSurfaces}
S.~Reuther, I.~Nitschke, A.~Voigt, A numerical approach for fluid deformable
  surfaces, Journal of Fluid Mechanics 900 (2020) R8.
\newblock \href {https://doi.org/10.1017/jfm.2020.564}
  {\path{doi:10.1017/jfm.2020.564}}.

\bibitem{olshanskii2018finite}
M.~A. Olshanskii, A.~Quaini, A.~Reusken, V.~Yushutin, A finite element method
  for the surface stokes problem, SIAM Journal on Scientific Computing 40
  (2018) A2492--A2518.

\bibitem{brandner2022finite}
P.~Brandner, T.~Jankuhn, S.~Praetorius, A.~Reusken, A.~Voigt, Finite element
  discretization methods for velocity-pressure and stream function formulations
  of surface stokes equations, SIAM Journal on Scientific Computing 44~(4)
  (2022) A1807--A1832.

\bibitem{CiCP-31-947}
M.~Nestler, A.~Voigt, Active nematodynamics on curved surfaces – the
  influence of geometric forces on motion patterns of topological defects,
  Communications in Computational Physics 31 (2022) 947--965.
\newblock \href {https://doi.org/10.4208/cicp.OA-2021-0206}
  {\path{doi:10.4208/cicp.OA-2021-0206}}.

\bibitem{WMB_arXiv_2022}
Z.~Wang, M.~C. Marchetti, F.~Brauns, Patterning of morphogenetic anisotropy
  fields, arXiv:2212.12215 (2022).

\bibitem{BKV_arXiv_2023}
E.~Bachini, V.~Krause, A.~Voigt, The interplay of geometry and coarsening in
  multicomponent lipid vesicles under the influence of hydrodynamics,
  arXiv:2302.04028 (2023).

\bibitem{grande2018analysis}
J.~Grande, C.~Lehrenfeld, A.~Reusken, Analysis of a high-order trace finite
  element method for {PDEs} on level set surfaces, SIAM Journal on Numerical
  Analysis 56 (2018) 228--255.

\bibitem{HarderingPraetorius2021TangentialErrors}
H.~Hardering, S.~Praetorius, Tangential errors of tensor surface finite
  elements, IMA Journal of Numerical Analysis (2022) drac015\href
  {https://doi.org/10.1093/imanum/drac015} {\path{doi:10.1093/imanum/drac015}}.

\bibitem{OlshanskiiEtAl2019InfSupStability}
M.~A. Olshanskii, A.~Reusken, A.~Zhiliakov, Inf-sup stability of the trace
  {P2-P1 Taylor-Hood} elements for surface {PDEs}, Mathematics of Computation
  90 (2021) 1527--1555.
\newblock \href {https://doi.org/10.1090/mcom/3551}
  {\path{doi:10.1090/mcom/3551}}.

\bibitem{JankuhnEtAl2021ErrorAnalysis}
T.~Jankuhn, M.~A. Olshanskii, A.~Reusken, A.~Zhiliakov, Error analysis of
  higher order trace finite element methods for the surface {S}tokes equation,
  Journal of Numerical Mathematics 29 (2021) 245--267.
\newblock \href {https://doi.org/10.1515/jnma-2020-0017}
  {\path{doi:10.1515/jnma-2020-0017}}.

\bibitem{bachini2022diffusion}
E.~Bachini, P.~Brandner, T.~Jankuhn, M.~Nestler, S.~Praetorius, A.~Reusken,
  A.~Voigt, Diffusion of tangential tensor fields: numerical issues and
  influence of geometric properties, arXiv:2205.12581 (2022).

\bibitem{ratz2006pde}
A.~R{\"a}tz, A.~Voigt, {PDE's} on surfaces---a diffuse interface approach,
  Communications in Mathematical Sciences 4 (2006) 575--590.

\bibitem{burger2009finite}
M.~Burger, Finite element approximation of elliptic partial differential
  equations on implicit surfaces, Computing and Visualization in Science 12
  (2009) 87--100.

\bibitem{hoffmann2022theory}
L.~A. Hoffmann, L.~N. Carenza, J.~Eckert, L.~Giomi, Theory of defect-mediated
  morphogenesis, Science Advances 8 (2022) eabk2712.

\bibitem{doi:10.1137/070708135}
A.~Demlow, Higher-order finite element methods and pointwise error estimates
  for elliptic problems on surfaces, SIAM Journal on Numerical Analysis 47
  (2009) 805--827.
\newblock \href {https://doi.org/10.1137/070708135}
  {\path{doi:10.1137/070708135}}.

\bibitem{cahn1958free}
J.~W. Cahn, J.~E. Hilliard, Free energy of a nonuniform system. {I.
  I}nterfacial free energy, The Journal of Chemical Physics 28~(2) (1958)
  258--267.

\bibitem{feng2003numerical}
X.~Feng, A.~Prohl, Numerical analysis of the {Allen-Cahn} equation and
  approximation for mean curvature flows, Numerische Mathematik 94~(1) (2003)
  33--65.

\bibitem{meshconv}
F.~Stenger, https://gitlab.mn.tu-dresden.de/iwr/meshconv.

\bibitem{ahrens200536}
J.~Ahrens, B.~Geveci, C.~Law, 36 - paraview: An end-user tool for large-data
  visualization, in: C.~D. Hansen, C.~R. Johnson (Eds.), Visualization
  Handbook, Butterworth-Heinemann, Burlington, 2005, pp. 717--731.
\newblock \href
  {https://doi.org/https://doi.org/10.1016/B978-012387582-2/50038-1}
  {\path{doi:https://doi.org/10.1016/B978-012387582-2/50038-1}}.

\bibitem{ayachit2015paraview}
U.~Ayachit, The paraview guide: a parallel visualization application, Kitware,
  Inc., 2015.

\bibitem{Sussmanetal_JCP_1994}
M.~Sussmann, P.~Smereka, S.~Osher, A level set approach for computing solutions
  to incompressible two-phase flow, J. Comput. Phys. 119 (1995) 146--159.

\bibitem{Russoetal_JCP_2000}
G.~Russo, P.~Smereka, A remark on computing distance functions, J. Comput.
  Phys. 163 (2000) 51--67.

\bibitem{Bornemannetal_CVS_2006}
F.~Bornemann, C.~Rasch, Finite-element discretization of static
  {Hamilton-Jacobi} equations based on a local variational principle, Comput.
  Vis. Sci. 9 (2006) 57--69.

\bibitem{lee2017revisiting}
B.~Lee, J.~Darbon, S.~Osher, M.~Kang, Revisiting the redistancing problem using
  the {Hopf--Lax} formula, Journal of Computational Physics 330 (2017)
  268--281.

\bibitem{royston2018parallel}
M.~Royston, A.~Pradhana, B.~Lee, Y.~T. Chow, W.~Yin, J.~Teran, S.~Osher,
  Parallel redistancing using the {Hopf--Lax} formula, Journal of Computational
  Physics 365 (2018) 7--17.

\bibitem{vey2007amdis}
S.~Vey, A.~Voigt, Amdis: adaptive multidimensional simulations, Computing and
  Visualization in Science 10 (2007) 57--67.

\bibitem{witkowski2015software}
T.~Witkowski, S.~Ling, S.~Praetorius, A.~Voigt, Software concepts and numerical
  algorithms for a scalable adaptive parallel finite element method, Advances
  in Computational Mathematics 41 (2015) 1145--1177.

\bibitem{petsc-user-ref}
S.~Balay, S.~Abhyankar, M.~F. Adams, S.~Benson, J.~Brown, P.~Brune,
  K.~Buschelman, E.~Constantinescu, L.~Dalcin, A.~Dener, V.~Eijkhout, W.~D.
  Gropp, V.~Hapla, T.~Isaac, P.~Jolivet, D.~Karpeev, D.~Kaushik, M.~G. Knepley,
  F.~Kong, S.~Kruger, D.~A. May, L.~C. McInnes, R.~T. Mills, L.~Mitchell,
  T.~Munson, J.~E. Roman, K.~Rupp, P.~Sanan, J.~Sarich, B.~F. Smith,
  S.~Zampini, H.~Zhang, H.~Zhang, J.~Zhang, {PETSc/TAO} users manual, Tech.
  Rep. ANL-21/39 - Revision 3.17, Argonne National Laboratory (2022).

\bibitem{burger2006surface}
M.~Burger, Surface diffusion including adatoms, Communications in Mathematical
  Sciences 4 (2006) 1--51.

\bibitem{ratz2006diffuse}
A.~R{\"a}tz, A.~Voigt, A diffuse-interface approximation for surface diffusion
  including adatoms, Nonlinearity 20 (2006) 177.

\bibitem{lowengrub2009phase}
J.~S. Lowengrub, A.~R{\"a}tz, A.~Voigt, Phase-field modeling of the dynamics of
  multicomponent vesicles: Spinodal decomposition, coarsening, budding, and
  fission, Physical Review E 79 (2009) 031926.

\bibitem{teigen2011diffuse}
K.~E. Teigen, P.~Song, J.~Lowengrub, A.~Voigt, A diffuse-interface method for
  two-phase flows with soluble surfactants, Journal of Computational Physics
  230 (2011) 375--393.

\bibitem{marth2014signaling}
W.~Marth, A.~Voigt, Signaling networks and cell motility: a computational
  approach using a phase field description, Journal of Mathematical Biology 69
  (2014) 91--112.

\bibitem{doi:10.1137/21M1433642}
P.~Werner, M.~Burger, F.~Frank, H.~Garcke, A diffuse interface model for cell
  blebbing including membrane-cortex coupling with linker dynamics, SIAM
  Journal on Applied Mathematics 82 (2022) 1091--1112.
\newblock \href {https://doi.org/10.1137/21M1433642}
  {\path{doi:10.1137/21M1433642}}.

\bibitem{aland2012buckling}
S.~Aland, A.~R{\"a}tz, M.~R{\"o}ger, A.~Voigt, Buckling instability of viral
  capsids—a continuum approach, Multiscale Modeling \& Simulation 10 (2012)
  82--110.

\bibitem{aland2011continuum}
S.~Aland, J.~Lowengrub, A.~Voigt, A continuum model of colloid-stabilized
  interfaces, Physics of Fluids 23~(6) (2011) 062103.

\bibitem{PhysRevE.86.046321}
S.~Aland, J.~Lowengrub, A.~Voigt, Particles at fluid-fluid interfaces: A new
  navier-stokes-cahn-hilliard surface- phase-field-crystal model, Phys. Rev. E
  86 (2012) 046321.
\newblock \href {https://doi.org/10.1103/PhysRevE.86.046321}
  {\path{doi:10.1103/PhysRevE.86.046321}}.

\bibitem{li2009solving}
X.~Li, J.~Lowengrub, A.~R{\"a}tz, Voigt, Solving {PDEs} in complex geometries:
  a diffuse domain approach, Communications in Mathematical Sciences 7 (2009)
  81--107.

\bibitem{teigen2009diffuse}
K.~E. Teigen, X.~Li, J.~Lowengrub, F.~Wang, A.~Voigt, A diffuse-interface
  approach for modeling transport, diffusion and adsorption/desorption of
  material quantities on a deformable interface, Communications in Mathematical
  Sciences 4 (2009) 1009--1037.

\bibitem{reuther2016incompressible}
S.~Reuther, A.~Voigt, Incompressible two-phase flows with an inextensible
  {N}ewtonian fluid interface, Journal of Computational Physics 322 (2016)
  850--858.

\end{thebibliography}
	
\end{document}